\documentclass{amsart}
\usepackage{latexsym}
\usepackage{amssymb}
\usepackage{amsmath}

\begin{document}


\newtheorem{theorem}{Theorem}
\newtheorem{problem}{Problem}
\newtheorem{definition}{Definition}
\newtheorem{lemma}{Lemma}
\newtheorem{proposition}{Proposition}
\newtheorem{corollary}{Corollary}
\newtheorem{example}{Example}
\newtheorem{conjecture}{Conjecture}
\newtheorem{algorithm}{Algorithm}
\newtheorem{exercise}{Exercise}
\newtheorem{remarkk}{Remark}

\newcommand{\be}{\begin{equation}}
\newcommand{\ee}{\end{equation}}
\newcommand{\bea}{\begin{eqnarray}}
\newcommand{\eea}{\end{eqnarray}}
\newcommand{\beq}[1]{\begin{equation}\label{#1}}
\newcommand{\eeq}{\end{equation}}
\newcommand{\beqn}[1]{\begin{eqnarray}\label{#1}}
\newcommand{\eeqn}{\end{eqnarray}}
\newcommand{\beaa}{\begin{eqnarray*}}
\newcommand{\eeaa}{\end{eqnarray*}}
\newcommand{\req}[1]{(\ref{#1})}

\newcommand{\lip}{\langle}
\newcommand{\rip}{\rangle}

\newcommand{\uu}{\underline}
\newcommand{\oo}{\overline}
\newcommand{\La}{\Lambda}
\newcommand{\la}{\lambda}
\newcommand{\eps}{\varepsilon}
\newcommand{\om}{\omega}
\newcommand{\Om}{\Omega}
\newcommand{\ga}{\gamma}
\newcommand{\rrr}{{\Bigr)}}
\newcommand{\qqq}{{\Bigl\|}}

\newcommand{\dint}{\displaystyle\int}
\newcommand{\dsum}{\displaystyle\sum}
\newcommand{\dfr}{\displaystyle\frac}
\newcommand{\bige}{\mbox{\Large\it e}}
\newcommand{\integers}{{\Bbb Z}}
\newcommand{\rationals}{{\Bbb Q}}
\newcommand{\reals}{{\rm I\!R}}
\newcommand{\realsd}{\reals^d}
\newcommand{\realsn}{\reals^n}
\newcommand{\NN}{{\rm I\!N}}
\newcommand{\DD}{{\rm I\!D}}
\newcommand{\degree}{{\scriptscriptstyle \circ }}
\newcommand{\dfn}{\stackrel{\triangle}{=}}
\def\complex{\mathop{\raise .45ex\hbox{${\bf\scriptstyle{|}}$}
     \kern -0.40em {\rm \textstyle{C}}}\nolimits}
\def\hilbert{\mathop{\raise .21ex\hbox{$\bigcirc$}}\kern -1.005em {\rm\textstyle{H}}} 
\newcommand{\RAISE}{{\:\raisebox{.6ex}{$\scriptstyle{>}$}\raisebox{-.3ex}
           {$\scriptstyle{\!\!\!\!\!<}\:$}}} 

\newcommand{\hh}{{\:\raisebox{1.8ex}{$\scriptstyle{\degree}$}\raisebox{.0ex}
           {$\textstyle{\!\!\!\! H}$}}}

\newcommand{\OO}{\won}
\newcommand{\calA}{{\mathcal A}}
\newcommand{\BB}{{\mathcal B}}
\newcommand{\calC}{{\cal C}}
\newcommand{\calD}{{\cal D}}
\newcommand{\calE}{{\cal E}}
\newcommand{\calF}{{\mathcal F}}
\newcommand{\calG}{{\cal G}}
\newcommand{\calH}{{\cal H}}
\newcommand{\calK}{{\cal K}}
\newcommand{\calL}{{\mathcal L}}
\newcommand{\calM}{{\cal M}}
\newcommand{\calO}{{\cal O}}
\newcommand{\calP}{{\cal P}}
\newcommand{\calU}{{\mathcal U}}
\newcommand{\calX}{{\cal X}}
\newcommand{\calXX}{{\cal X\mbox{\raisebox{.3ex}{$\!\!\!\!\!-$}}}}
\newcommand{\calXXX}{{\cal X\!\!\!\!\!-}}
\newcommand{\gi}{{\raisebox{.0ex}{$\scriptscriptstyle{\cal X}$}
\raisebox{.1ex} {$\scriptstyle{\!\!\!\!-}\:$}}}
\newcommand{\intsim}{\int_0^1\!\!\!\!\!\!\!\!\!\sim}
\newcommand{\intsimt}{\int_0^t\!\!\!\!\!\!\!\!\!\sim}
\newcommand{\pp}{{\partial}}
\newcommand{\al}{{\alpha}}
\newcommand{\sB}{{\cal B}}
\newcommand{\sL}{{\cal L}}
\newcommand{\sF}{{\cal F}}
\newcommand{\sE}{{\cal E}}
\newcommand{\sX}{{\cal X}}
\newcommand{\R}{{\rm I\!R}}
\renewcommand{\L}{{\rm I\!L}}
\newcommand{\vp}{\varphi}
\newcommand{\N}{{\rm I\!N}}
\def\ooo{\lip}
\def\ccc{\rip}
\newcommand{\ot}{\hat\otimes}
\newcommand{\rP}{{\Bbb P}}
\newcommand{\bfcdot}{{\mbox{\boldmath$\cdot$}}}

\renewcommand{\varrho}{{\ell}}
\newcommand{\dett}{{\textstyle{\det_2}}}
\newcommand{\sign}{{\mbox{\rm sign}}}
\newcommand{\TE}{{\rm TE}}
\newcommand{\TA}{{\rm TA}}
\newcommand{\E}{{\rm E\,}}
\newcommand{\won}{{\mbox{\bf 1}}}
\newcommand{\Lebn}{{\rm Leb}_n}
\newcommand{\Prob}{{\rm Prob\,}}
\newcommand{\sinc}{{\rm sinc\,}}
\newcommand{\ctg}{{\rm ctg\,}}
\newcommand{\loc}{{\rm loc}}
\newcommand{\trace}{{\,\,\rm trace\,\,}}
\newcommand{\Dom}{{\rm Dom}}
\newcommand{\ifff}{\mbox{\ if and only if\ }}
\newcommand{\nproof}{\noindent {\bf Proof:\ }}
\newcommand{\remark}{\noindent {\bf Remark:\ }}
\newcommand{\remarks}{\noindent {\bf Remarks:\ }}
\newcommand{\note}{\noindent {\bf Note:\ }}

\newcommand{\boldx}{{\bf x}}
\newcommand{\boldX}{{\bf X}}
\newcommand{\boldy}{{\bf y}}
\newcommand{\boldR}{{\bf R}}
\newcommand{\uux}{\uu{x}}
\newcommand{\uuY}{\uu{Y}}

\newcommand{\limn}{\lim_{n \rightarrow \infty}}
\newcommand{\limN}{\lim_{N \rightarrow \infty}}
\newcommand{\limr}{\lim_{r \rightarrow \infty}}
\newcommand{\limd}{\lim_{\delta \rightarrow \infty}}
\newcommand{\limM}{\lim_{M \rightarrow \infty}}
\newcommand{\limsupn}{\limsup_{n \rightarrow \infty}}

\newcommand{\ra}{ \rightarrow }

\newcommand{\ARROW}[1]
  {\begin{array}[t]{c}  \longrightarrow \\[-0.2cm] \textstyle{#1} \end{array} }

\newcommand{\AR}
 {\begin{array}[t]{c}
  \longrightarrow \\[-0.3cm]
  \scriptstyle {n\rightarrow \infty}
  \end{array}}

\newcommand{\pile}[2]
  {\left( \begin{array}{c}  {#1}\\[-0.2cm] {#2} \end{array} \right) }

\newcommand{\floor}[1]{\left\lfloor #1 \right\rfloor}

\newcommand{\mmbox}[1]{\mbox{\scriptsize{#1}}}

\newcommand{\ffrac}[2]
  {\left( \frac{#1}{#2} \right)}

\newcommand{\one}{\frac{1}{n}\:}
\newcommand{\half}{\frac{1}{2}\:}

\def\le{\leq}
\def\ge{\geq}
\def\lt{<}
\def\gt{>}

\def\squarebox#1{\hbox to #1{\hfill\vbox to #1{\vfill}}}
\newcommand{\nqed}{\hspace*{\fill}
           \vbox{\hrule\hbox{\vrule\squarebox{.667em}\vrule}\hrule}\bigskip}

\newcommand{\no}{\noindent}
\newcommand{\EE}{\mathbb{E}}
\newcommand{\RR}{\mathbb{R}}
\newcommand{\D}{\mathcal{D}}
\newcommand{\LL}{\mathcal{L}}
\newcommand{\F}{\mathcal{F}}
\newcommand{\B}{\mathcal{B}}
\newcommand{\M}{\mathcal{M}}
\newcommand{\W}{\mathbb{W}}
\newcommand{\grandX}{\mathbb{X}}
\newcommand{\loigrandX}{\mu^{\mathbb{X}}}

\title{VARIATIONAL CALCULUS FOR DIFFUSIONS}

\author{ K\'evin HARTMANN}
\maketitle
\noindent
{\bf Abstract:}{\small{We expand the classic variational formulation of $-\log\EE\left[e^{-f}\right]$ to the case where f depends on a diffusion, and not only a on Brownian motion, while decreasing the integrability hypothesis on f. We also give an entropic characterisation of the invertibility of a perturbation of a diffusion and discuss the attainability of the infimum in the aforementioned variational formulation.
}}\\

\vspace{0.5cm}
\no
\tableofcontents

\section{\bf{Introduction}}

Denote $\W$ the space of continuous functions from $[0,1]$ to $\RR^n$ and H the associated canonical Cameron-Martin space of elements of $\W$ which admit a density in $L^2$. Also denote $\mu$ the Wiener measure, W the coordinate process, and $(\F_t)$ the canonical filtration of W completed with respect to $\mu$. W is a Brownian motion under $\mu$. Set f a bounded from above measurable function from W to $\RR$. In \cite{du}, Dupuis and Ellis prove that
\bea -\log\EE_\mu\left[e^{-f}\right]=\inf_\theta\left(\EE_\theta\left[f\right]+H(\theta|\mu)\right)\label{1rv0}\eea
\no where the infimum is taken over the probability measures $\theta$ on $\W$ which are absolutely continuous with respect to $\mu$ and the relative entropy $H(\theta|\mu)$ is equal to $\EE_\mu\left[\frac{d\theta}{d\mu}\log\frac{d\theta}{d\mu}\right]$. In \cite{bd}, Bou\'e and Dupuis use it to derive the variational formulation
\bea -\log\EE_\mu\left[e^{-f}\right]=\inf_u\EE_\mu\left[f\circ (W+u)+\frac{1}{2}\int_0^1|\dot{u}(s)|^2ds\right]\label{1rv}\eea
\no where the infimum is taken over $L^2$ functions from $\W$ to H whose density is adapted to $(\F_t)$. This variational formulation is useful to derive large deviation asymptotics as Laplace principles for small noise diffusions for instance. This result was later extended by Budhiraja and Dupuis to Hilbert-space-valued Brownian motions in \cite{bud}, and then by Zhang to abstract Wiener spaces in \cite{zh}, using the framework developed by \"Ust\"unel and Zakai in \cite{ust1}.\newline
\no The bounded from above hypothesis in \ref{1rv} was weakened significantly by \"Ust\"unel in \cite{art}, it was replaced with the condition
$$\EE_\mu\left[fe^{-f}\right]<\infty$$
\no and the existence of conjugate integers p and q such that
$$f\in L^p(\mu),e^{-f}\in L^q(\mu)$$
\no These relaxed hypothesis pave the way to new applications. The possibility of using unbounded functions is primordial in Dabrowski's application of \ref{1rv} to free entropy in \cite{da}.\newline
\no \"Ust\"unel's approach is routed in the study of the perturbations of the identity of $\W$, which is the coordinate process, and their invertibility. The question of the invertibility of an adapted perturbation of the identity is linked to the representability of measures and was put to light by the celebrated example of Tsirelson \cite{tsi}. \"Ust\"unel proved that if $u\in L^2(\mu,H)$ and has an adapted density, $I_\W+u$ is $\mu$-a.s. invertible if and only if
$$H((I_\W+u)\mu|\mu)=\frac{1}{2}\EE_\mu\left[|u|_H^2\right]$$
\no To prove \ref{1rv} with the integrability conditions specified above, \"Ust\"unel uses the fact that H-$C^1$ shifts, meaning shifts that are a.s. Fr\'echet-differentiable on H with an a.s. continuous on H Fr\'echet derivative, are a.s. invertible, and that shifts can be approached with H-$C^1$ shifts using the Ornstein-Uhlenbeck semigroup.\newline
This paper focuses on getting a variational formulation similar as the one above in the case of a diffusion V which satisfies a stochastic differential equation
$$V(t)=c+\int_0^t\sigma(V(s))dB(s)+\int_0^tb(V(s))ds$$
\no where B is a Brownian motion, thus generalizing the case of the Brownian motion. We also weaken the integration hypothesis on f since we only require $\EE[fe^{-f}]<\infty$ and $f\in L^p(\mu)$ for some $p>1$. \"Ust\"unel's proof consists in approaching f with H-$C^1$ functions and then use H-$C^1$ shifts, which are invertible, obtained using those functions. This approach is deeply rooted in the Brownian motion specific case, since it relies on sophisticated stochastic analysis tool that were developed for a Gaussian framework. Here we write the density  $\frac{e^{-f}}{\EE[e^{-f}]}$ as the Wick exponential of some v and then approach v with retarded shifts which generate invertible perturbations of the identity. Since we work under the law of a diffusion and not the Wiener measure, the perturbations of the identity we consider are not affine shifts. We work on $\W$ under the image measure of $(V,B)$ that we denote $\loigrandX$ and we construct a Brownian motion $\beta_\grandX$ such that W verifies
$$W(t)=c+\int_0^t\sigma(W(s))d\beta_\grandX(s)+\int_0^tb(W(s))ds$$
\no We only consider perturbations that verify the Girsanov condition. If u is such a perturbation, we denote $X^u$ the solution of the stochastic differential equation
$$X^u=c+\int_0^t\sigma(X^u(s))d\left(\beta_\grandX+u\right)(s)+\int_0^tb(X^u(s))ds$$
\no and $\grandX^u=(X^u,\beta_\grandX+u)$. $\grandX^u$ plays the same role as $W+u$ in the Brownian case and it is invertible if and only if
$$H(\grandX^u\loigrandX|\loigrandX)=\frac{1}{2}\EE_{\loigrandX}\left[|u|_H^2\right]$$
\no We conclude the paper with a discussion over the attainability of the infimum in the variational formulation.\newline

\section{\bf{Framework}}

\no Set $m\leq d\in \NN^*$, $c\in\RR^m$, $\sigma:\RR^m\rightarrow\M_{m,d}(\RR)$ and $b:\RR^m\rightarrow\RR^m$ bounded and lipschitz functions. $\sigma_i$ will denote the
i-th column of $\sigma$. Notice that every matrix
will be identified with its canonical linear operator. Set $(\Omega,R,(\mathcal{G}_t))$ a probability space, V a R-Brownian motion on
$\Omega$ with values in $\RR^d$. Set Y a $\RR^m$-valued strong solution of the
stochastic differential equation:
$$Y(t)=c+\int_0^t\sigma(Y(s))dV(s)+\int_0^tb(Y(s))ds$$
\no on $(\Omega,R,(\mathcal{G}_t),B)$. The hypotheses on $\sigma$ and b ensure the
existence and uniqueness of Y if we impose its paths to be continuous.\newline
\no We denote $\mu$
the Wiener measure on $C([0,1],\RR^d)$ and $\mu^X$ the image measure of X.  We denote $\W=C([0,1],\RR^{m+d})$ and we consider the measure $\mu^X\times\mu$ on $\W$.\newline
\no We define the processes X and B on $\W$ by:
\beaa X(t)&:&(w,w')\in \W\mapsto w(t)\in\RR^m\\
B(t)&:&(w,w')\in \W\mapsto w'(t)\in\RR^d\eeaa
\no Under $\mu^X\times\mu$, the law of X is $\mu^X$, B is a Brownian
motion and they are independent. We denote $X_i$ and $B_i$ the i-th
coordinates of x and B.
\no Define $M=X-c-\int_0^.b(X(s))ds$ and $a=\sigma\sigma^T$. For $1\leq i\leq d$,
$M^i=X^i-c_i-\int_0^.b^i(X(s))ds$ is a local martingale and we have:
$$\langle M^i,M^j\rangle=\int_0^.a^{ij}(X(s))ds$$
\no Now set $y\in\RR$. observe that
\beaa \tilde{\sigma}(y):\begin{array}{rcl}
  ker(\sigma(y))^\perp&\rightarrow& im(\sigma(y))\\
e&\mapsto&\sigma(y)(e)\end{array}\eeaa
\no is an isomorphism and set $\theta(y)$ the unique element of $\mathcal{M}_{d,m}(\RR)$
which is equal to $\tilde{\sigma}(y)^{-1}$ on $im(\sigma(y))$ and 0 on
$im(\sigma(y))^\perp$ and $\eta(y)$ the unique element of
$\mathcal{M}_d(\RR)$ which is equal to 0 on $\ker(\sigma(y))^\perp$
and to the identity on $\ker(\sigma(y))$.\newline
\no Notice that we have $(\theta\sigma+\eta)(y)=I_d(\RR)$.\newline
\no We define
$$\beta_\grandX=\int_0^.\theta(X(s))dM(s)+\int_0^.\eta(X(s))dB(s)$$
\no $\beta_{\grandX,i}$ will denote the i-th coordinate of $\beta_\grandX$.\newline
\no $\beta_\grandX$ is a Brownian motion. Indeed, it is clearly a local
martingale and since m and B are independent:
\beaa \left(\langle
  \beta_{\grandX,i},\beta_{\grandX,j}\rangle(t)\right)_{i,j}&=&\int_0^t\left(\theta(X(s))\;\;\eta(X(s))\right)\left(\begin{array}{rcl}
    a(X(s))&\;&0\\0&\;&I\end{array}\right)\left(\theta(X(s))\;\;\eta(X(s))\right)^Tds\\
&=&\int_0^t\theta(X(s))\sigma(X(s))\sigma(X(s))^T\theta(X(s))^T+\eta(X(s))\eta(X(s))^Tds\\
&=&tI_d(\RR)\eeaa
\no Moreover $M=\int_0^.\sigma(X(s))d\beta_\grandX(s)$. Indeed,
$M-\int_0^.\sigma(X(s))d\beta_\grandX(s)$ is a local martingale and:
\beaa &&\left(\left\langle
  M-\int_0^.\sigma(X(s))d\beta_\grandX(s),M-\int_0^.\sigma(X(s))d\beta_\grandX(s)\right\rangle(t)\right)_{i,j}\\
&=&\int_0^t\left(\sigma(X(s))\theta(X(s))+I\;\;\sigma(X(s))\eta(X(s))\right)\\
&&\left(\begin{array}{rcl}
    a(X(s))&\;&0\\0&\;&I\end{array}\right)\left(\sigma(X(s))\theta(X(s))+I\;\;\sigma(X(s))\eta(X(s))\right)^Tds\\
&=&\int_0^t\left(\sigma(X(s))\theta(X(s))+I\right)\sigma(X(s))\sigma(X(s))^T\left(\sigma(X(s))\theta(X(s))-I\right)^T\\
&&+\sigma(X(s))\eta(X(s))\eta(X(s))^T\eta(X(s))^Tds\\
&=&0\eeaa
\no This construction of $\beta_\grandX$ is taken from \cite{ro}.\newline
\no We denote
$$\grandX=(X,\beta_\grandX)$$
\no and $\loigrandX$ its image measure.
\no X is a $\loigrandX$ path-continuous strong solution of the
stochastic differential equation
\beaa X=&c+\int_0^.\sigma(X(s))d\beta_\grandX(s)+\int_0^.b(X(s))ds\eeaa
\no The filtration of a process m will be denoted
$\left(\F^m_t\right)$, the filtration of $\grandX$ will be simply denoted
$\left(\F_t\right)$. Except if stated otherwise, every filtration
considered is completed with respect to $\loigrandX$.
\no If m is a martingale and v admits a density $\dot{v}$ whose stochastic
integral with respect to m is well defined we will denote
$$\delta_m v =\int_0^1\dot{v}(s)dm(s)$$
\no We also denote the Wick exponential as follow
$$\rho(\delta_m
v)=\exp\left(\int_0^1\dot{v}(s)dm(s)-\frac{1}{2}\int_0^1\left|\dot{v}(s)\right|^2d\langle
  m\rangle(s)\right)$$
\no and for $p\geq 0$ we denote
\beaa L^p_a(\loigrandX,H)&=&\left\{u\in L^p(\loigrandX,H), \dot{u}\;is\;(\F_t)-adapted\right\}\\
\\G_p(\loigrandX,m)&=&\left\{u\in L_a^p(\loigrandX,H), \EE_{\loigrandX}\left[\rho(-\delta_m u)\right]=1\right\}\eeaa
\no We denote $H=\left\{\int_0^.\dot{h}(s)ds, \dot{h}\in
  L^2([0,1],\RR^d)\right\}$. For $u\in G_0(\loigrandX,\beta_\grandX)$, we define
$\beta_\grandX^u:=\beta_\grandX+u$ and $X^u$ a path-continuous strong solution
of the stochastic differential equation
$$X^u=c+\int_0^.\sigma(X^u(s))d\beta_\grandX^u(s)+\int_0^.b(X^u(s))ds$$
\no on $(W,\loigrandX,(\F_t),\beta_\grandX)$. Once again the hypotheses on $\sigma$
ensure the existence and $\loigrandX$-path uniqueness of $X^u$. We also denote
$$M^u=X^u-c-\int_0^.b(X^u(s))ds=\int_0^.\sigma(X^u(s))d\beta_\grandX^u(s)$$
\no and
\beaa
\grandX^u&=&(X^u,\beta_\grandX+u)\eeaa

\no We have a Girsanov-like change of measure theorem relative to
$\loigrandX$:

\begin{proposition}
\label{1gir}
Set $u\in G_0(\loigrandX,\beta_\grandX)$ , for every bounded Borel function f:
$$\EE_{\loigrandX}\left[f\right]=\EE_{\loigrandX}\left[f\circ \grandX^u\rho(-\delta_{\beta_\grandX} u)\right]$$
\end{proposition}
\nproof Set f a bounded Borel function and $u\in G_0(\loigrandX,\beta_\grandX)$, denote
$\theta$  the probability on $\W$ defined by
 $$\frac{d\theta}{d\loigrandX}=\rho(-\delta_{\beta_\grandX}
  u)$$
\no According to the Girsanov theorem, the law of
  $\beta_\grandX+u$ under $\theta$ is the same as the law of
  $\beta_\grandX$ under $\loigrandX$. Consequently, the law of $X^u$ under $\theta$ is the same as
  the law of $X$ under $\loigrandX$ and
$$\EE_{\loigrandX}\left[f\circ \grandX\right]=\EE_{\loigrandX}\left[f\circ \grandX^u\rho(-\delta_{\beta_\grandX} u)\right]$$\nqed

\begin{theorem}
\no Set $u\in G_0(\loigrandX,\beta_\grandX)$, we have
$$\grandX^u\loigrandX\sim \loigrandX$$
\end{theorem}

\nproof Set $f\in C_b(\W)$ and set $\theta$ the measure on $\W$ given by
$$\frac{d\theta}{d\loigrandX}=\rho(-\delta_{\beta_\grandX} u)$$
\no We have
\beaa\EE_{\grandX^u\theta}\left[f\right]&=&\EE_\theta\left[f\circ \grandX^u\right]\\
&=&\EE_{\loigrandX}\left[f\circ\grandX^u\rho(-\delta_{\beta_\grandX} u)\right]\\
&=&\EE_{\loigrandX}\left[f\right]\eeaa
\no so $\grandX^u\theta=\loigrandX$.\newline\no Since $\theta\sim\loigrandX$, $\grandX^u\theta\sim\grandX^u\loigrandX$, which conclude the proof.\nqed

\no Set $u,v\in G_0(\loigrandX,\beta_\grandX)$, this theorem ensures that
if g is a random variable defined on $\W$,  the composition $g\circ \grandX^v$ is
well-defined. Indeed set $\tilde{g}$ and $\hat{g}$ in the same
equivalence class in $L^0(\loigrandX)$. Then
\beaa \loigrandX\left(\hat{g}\circ \grandX^v=\tilde{g}\circ
  \grandX^v\right)=\grandX^v\loigrandX\left(\tilde{g}=\hat{g}\right)=1\eeaa
\no since $\grandX^v\loigrandX\ll \loigrandX$.\newline
\no In particular the compositions $u\circ \grandX^v$ and $\grandX^u\circ \grandX^v$ are well-defined
since u and $\grandX^u$ are random variables defined on $W$ with values in $H$
and $W$ respectively.\newline

\section{\bf{Action of the composition by $\grandX^u$ and invertibility results}}

\begin{proposition}
\label{1leminv}
\no Set $u\in G_0(\loigrandX,\beta_\grandX)$. We have $\loigrandX$-a.s.:
\beaa M^u&=&M\circ \grandX^u\\
\beta_\grandX^u&=&\beta_\grandX\circ \grandX^u\eeaa
\end{proposition}

\nproof
We have
\beaa M\circ \grandX^u&=&\left(X-c-\int_0^.b(X(s))ds\right)\circ \grandX^u\\
&=&X^u-c-\int_0^.b(X^u(s))ds\\
&=&M^u\eeaa
Now,

\beaa \beta_\grandX\circ \grandX^u
&=&\left(\int_0^.\theta(X(s))dM(s)+\int_0^.\eta(X(s))dB(s)\right)\circ \grandX^u\\
&=&\int_0^.\theta(X^u(s))dM^u(s)+\int_0^.\eta(X^u(s))d\beta_\grandX^u(s)\\
&=&\int_0^.\theta(X^u(s))\sigma(X^u(s))d\beta_\grandX^u(s)+\int_0^.\eta(X^u(s))d\beta_\grandX^u(s)\\
&=&\beta_\grandX^u\eeaa

\begin{proposition}
\no Set $u,v\in G_0(\loigrandX,\beta_\grandX)$ such that $v+u\circ\grandX^u\in G_0(\loigrandX,\beta_\grandX)$, we have $\loigrandX$-a.s.:
$$\grandX^u\circ \grandX^v=\grandX^{v+u\circ \grandX^v}$$
\end{proposition}

\no We have
\beaa \grandX^u\circ \grandX^v = (X^u\circ \grandX^v,(\beta_\grandX+u)\circ \grandX^v)=(X^u\circ \grandX^v,
\beta_\grandX+v+u\circ \grandX^u)\eeaa
\no Now,
\beaa X^u\circ \grandX^v &=&
\left(c+\int_0^.\sigma(X^u(s))d\beta_\grandX^u(s)+\int_0^.b(X^u(s))ds\right)\circ \grandX^v\\
&=&\left(c+\int_0^.\sigma(X^u(s))d\beta_\grandX(s)+\int_0^.\sigma(X^u(s))du(s)+\int_0^.b(X^u(s))ds\right)\circ
\grandX^v\\
&=&c+\int_0^.\sigma(X^u(s)\circ \grandX^v)d\beta_\grandX^v(s)+\int_0^.\sigma(X^u(s)\circ
\grandX^v)\dot{u}(s)\circ \grandX^v ds+\int_0^.b(X^u(s)\circ \grandX^v)ds\\
&=&c+\int_0^.\sigma(X^u(s)\circ \grandX^v)d\beta_\grandX^{v+u\circ
  \grandX^v}(s)+\int_0^.b(X^u(s)\circ \grandX^v)ds\eeaa
\no $X^u\circ \grandX^v$ and $X^{v+u\circ \grandX^v}$ are path continuous strong
solutions to the same stochastic differential equation so they are
equal $\loigrandX$-a.s.\newline
\no Finally, we have $\loigrandX$-a.s.
\beaa \grandX^u\circ \grandX^v=(X^{v+u\circ \grandX^v},\beta_\grandX+v+u\circ \grandX^v)=\grandX^{v+u\circ \grandX^v}\eeaa\nqed

\section{\bf{Invertibility results}}

\begin{definition}
A measurable map $U:\W\rightarrow \W$ is said to be
$\loigrandX$-a.s. left-invertible if and only if $U\loigrandX\ll\loigrandX$ and there exists a measurable map
$V:\W\rightarrow \W$ such that $V\circ U=I_\W$ $\loigrandX$-a.s.\newline
\no A measurable map $U:\W\rightarrow \W$ is said to be
$\loigrandX$-a.s. right-invertible if and only if there exists a measurable map
$V:\W\rightarrow \W$ such that $V\loigrandX\ll\loigrandX$ and  $U\circ V=I_\W$ $\loigrandX$-a.s.
\end{definition}

\begin{proposition}
\label{1lrinv}
Set $U,V:\W\rightarrow\W$ measurable maps such that $V\circ U=I_\W$
$\loigrandX$-a.s. and $V\loigrandX\ll\loigrandX$ Then  $U\circ V=I_\W$ $U\loigrandX$-a.s., so if $U\loigrandX\sim \loigrandX$, we
also have $U\circ V=I_\W$ $\loigrandX$-a.s. In that case, we will say that
$U$ is $\loigrandX$-a.s. invertible.
\end{proposition}

\nproof There exists $A\subset W$ such that $\loigrandX(A)=1$ and for every
$w\in A$, $V\circ U(w)=w$. Consider such a set A, we have
\beaa \EE_{U\loigrandX}\left[1_{U\circ
    V(w)=w}\right]&=&\EE_{\loigrandX}\left[1_{U\circ V\circ U(w)=U(w)}\right]\\
&=&\EE_{\loigrandX}\left[1_{U\circ V\circ U(w)=U(w)}1_{w\in
      A}\right]+\EE_{\loigrandX}\left[1_{U\circ V\circ U(w)=U(w)}1_{w\notin
        A}\right]\\
&=&\EE_{\loigrandX}\left[1_{ U(w)=U(w)}1_{w\in
      A}\right]\\
&=&1\eeaa\nqed

\section{\bf{Entropic characterisation of the invertibility of $\grandX^u$}}

In this section, we prove that the process $\grandX^u$ is left invertible if
and only if the kinetic energy of the perturbation u is equal to the relative entropy of $\grandX^u\loigrandX$.

\begin{proposition}
\no Set $u\in G_2(\loigrandX,\beta_\grandX)$. We have:
$$H(\grandX^u\loigrandX|\loigrandX)\leq\frac{1}{2}\EE_{\loigrandX}\left[|u|^2_H\right]$$
\end{proposition}

\nproof
Set $g\in C_b(\W)$ and denote $L=\frac{d\grandX^u\loigrandX}{d\loigrandX}$, we
have:
\beaa \EE_{\loigrandX}\left[g\circ \grandX^u\right]&=&\EE_{\loigrandX}\left[gL\right]\\
&=&\EE_{\loigrandX}\left[g\circ \grandX^u L\circ \grandX^u\rho(-\delta_{\beta_\grandX} u)\right]\eeaa
\no So $\loigrandX$-a.s.
$$L\circ \grandX^u\EE_{\loigrandX}\left[\left.\rho(-\delta_{\beta_\grandX} u)\right|\F_1^{\grandX^u}\right]= 1$$
\no and
\beaa H(\grandX^u\loigrandX|\loigrandX)&=&\EE_{\loigrandX}\left[L\log L\right]\\
&=&\EE_{\grandX^u\loigrandX}\left[\log L\right]\\
&=&\EE_{\loigrandX}\left[\log L\circ \grandX^u\right]\\
&=&-\EE_{\loigrandX}\left[\log\EE_{\loigrandX}\left[\left.\rho(-\delta_{\beta_\grandX} u)\right|\F_1^{\grandX^u}\right]\right]\\
&\leq&-\EE_{\loigrandX}\left[\log\rho(-\delta_{\beta_\grandX} u)\right]\\
&\leq&\frac{1}{2}\EE_{\loigrandX}\left[|u|_H^2\right]\eeaa
\nqed

\no Now comes the criteria:

\begin{theorem}
\label{1inv}
Set $u\in G_2(\loigrandX,\beta_\grandX)$. The three following propositions are equivalent:\newline
\no (i) $H(\grandX^u\loigrandX|\loigrandX)=\frac{1}{2}\EE_{\loigrandX}\left[|u|_H^2\right]$\newline
(ii) There exists $v\in G_0(\loigrandX,\beta_\grandX)$ such that $\grandX^v\circ \grandX^u=\grandX^u\circ \grandX^v=I_\W$ $\loigrandX$-a.s.\newline
\no (iii) $\grandX^u$ is $\loigrandX$-a.s. left-invertible
\end{theorem}

\nproof We first prove $(i)\Rightarrow (ii)$. We still denote
$L=\frac{d\grandX^u\loigrandX}{d\loigrandX}$ and as in the proof of last proposition we have $\loigrandX$-a.s.
$$L\circ \grandX^u\EE_{\loigrandX}\left[\rho(-\delta_{\beta_\grandX} u)|\grandX^u\right]= 1$$

\no Using Jensen inequality we have $\loigrandX$-a.s.

\beaa 0&=& \log L\circ
\grandX^u+\log\EE_{\loigrandX}\left[\left.\rho(-\delta_{\beta_\grandX}
    u)\right|\F^{\grandX^u}_1\right]\\
&\geq& \log L\circ
\grandX^u+\EE_{\loigrandX}\left[\left. \log\rho(-\delta_{\beta_\grandX}
    u)\right|\F^{\grandX^u}_1\right]\eeaa
\no and
\beaa 0&\geq&\EE_{\loigrandX}\left[\log L\circ \grandX^u\right]+\EE_{\loigrandX}\left[\log\rho(-\delta_{\beta_\grandX}
    u)\right]\\
&\geq&H(\grandX^u\loigrandX|\loigrandX)-\frac{1}{2}\EE_{\loigrandX}\left[|u|_H^2\right]\\
&=&0\eeaa
\no So
\beaa 0&=&\log L\circ \grandX^u+\log\EE_{\loigrandX}\left[\left.\rho(-\delta_{\beta_\grandX} u)\right|\F_1^{\grandX^u}\right]\\
&=&\log L\circ \grandX^u+\EE_{\loigrandX}\left[\log\left.\rho(-\delta_{\beta_\grandX} u)\right|\F_1^{\grandX^u}\right]\eeaa
\no and
$$\log\EE_{\loigrandX}\left[\left.\rho(-\delta_{\beta_\grandX} u)\right|\F_1^{\grandX^u}\right]=\EE_{\loigrandX}\left[\log\left.\rho(-\delta_{\beta_\grandX} u)\right|\F_1^{\grandX^u}\right]$$
\no The strict concavity of the function $\log$ gives
$$\EE_{\loigrandX}\left[\left.\rho(-\delta_{\beta_\grandX} u)\right|\F_1^{\grandX^u}\right]=\rho(-\delta_{\beta_\grandX} u)$$
\no Finally we have
\begin{equation}\label{1eg}L\circ \grandX^u\rho(-\delta_{\beta_\grandX} u)= 1\end{equation}

\no Since $\beta_\grandX$ is a $\loigrandX$-Brownian motion, there exists $v\in
G_0(\loigrandX,\beta_\grandX)$ such that
$L=\rho(-\delta_{\beta_\grandX} v)$.
\newline
\no We apply the logarithm to \ref{1eg} to get:
\beaa 0&=&\delta_{\beta_\grandX} v\circ \grandX^u+\frac{1}{2}|v\circ \grandX^u|^2_H+\delta_{\beta_\grandX}
u+\frac{1}{2}|u|^2_H\eeaa
\no We have:
\beaa \delta_{\beta_\grandX} v\circ \grandX^u &=&\int_0^1\dot{v}(s)\circ \grandX^ud\beta_\grandX(s) +\langle v\circ \grandX^u,u\rangle_H\eeaa

\no so finally we have:
\begin{equation}\label{1btest} 0= \delta_{\beta_\grandX}(v\circ \grandX^u+u)+\frac{1}{2}|v\circ \grandX^u+u|_H^2\end{equation}
According to Girsanov theorem $\beta_\grandX+v$ is a $\grandX^u\loigrandX$-Brownian motion,
so:
\beaa \EE_{\loigrandX}\left[L\log L\right] &=& \EE_{\grandX^u\loigrandX}\left[\log
    L\right]\\
&=&\EE_{\grandX^u\loigrandX}\left[-\int_0^1\dot{v}(s)d\beta_\grandX(s)-\frac{1}{2}\int_0^1|\dot{v}(s)|^2ds\right]\\
&=&\frac{1}{2}\EE_{\grandX^u\loigrandX}\left[\int_0^1|\dot{v}(s)|^2ds\right]\\
&=&\frac{1}{2}\EE_{\loigrandX}\left[|v\circ \grandX^u|_H^2\right]\eeaa

\no So $v\circ \grandX^u\in L^2_a(\loigrandX,H)$ and we can take the expectation
with respect to $\nu$ in \ref{1btest} to obtain $u+v\circ \grandX^u=0$
$\loigrandX$-a.s. and $\grandX^v\circ \grandX^u=I_\W$
$\loigrandX$-a.s. and $\grandX^v$ is a left-inverse of $\grandX^u$.\newline
\no Since $\grandX^v\loigrandX\sim\loigrandX$, we also have $\grandX^u\circ\grandX^v=I_W$ $\loigrandX$-a.s. from proposition \ref{1lrinv}.\newline
\no $(ii)\Rightarrow (iii)$ is immediate. Now we prove $(iii)\Rightarrow (i)$. We still denote $L=\frac{d\grandX^u\loigrandX}{d\loigrandX}$.\newline

\no Assume that $\grandX^u$ admits a left inverse V. Set
$v=-u\circ V$.\newline
\no We have $\loigrandX$-a.s.
$$v\circ\grandX^u=-u$$
\no and
$$\EE_{\grandX^u\loigrandX}\left[1_{\int_0^1|\dot{v}(s)|^2ds<\infty}\right]=\EE_{\loigrandX}\left[1_{\int_0^1|\dot{u}(s)|^2ds<\infty}\right]=1$$
\no so $v\in L^0(\grandX^u\loigrandX,H)$ and $v\in L^0(\loigrandX,H)$ since $\grandX^u\loigrandX\sim\loigrandX$.\newline
\no Now set $\dot{v}^n=\max(n,\min(\dot{v},-n))$, $\dot{v}^n\circ
\grandX^u$ is adapted. Set $A\in L^2(dt\times d\loigrandX)$ an adapted process, we
have:
\beaa\EE_{\loigrandX}\left[\rho(-\delta_{\beta_\grandX} u)\int_0^1\dot{v}^n(s)\circ \grandX^u
  A(s)\circ \grandX^uds\right]&=&\EE_{\loigrandX}\left[\int_0^1\dot{v}^n(s)
  A(s)ds\right]\\
&=&\EE_{\loigrandX}\left[\int_0^1\EE_{\loigrandX}\left[\left.\dot{v}^n(s)\right|\F^\grandX(s)\right]
  A(s)ds\right]\\
&=&\EE_{\loigrandX}\left[\rho(-\delta_{\beta_\grandX} u)\int_0^1\EE_{\loigrandX}\left[\left.\dot{v}^n(s)\right|\F^\grandX(s)\right]\circ \grandX^u
  A(s)\circ \grandX^uds\right]\eeaa
\no So $\EE_{\loigrandX}\left[\left.\dot{v}^n(s)\right|\F^\grandX(s)\right]\circ \grandX^u=\dot{v}^n(s)\circ
\grandX^u$ $dt\times d\loigrandX$-a.s. which implies
$\EE_{\loigrandX}\left[\left.\dot{v}^n(s)\right|\F^\grandX(s)\right]=\dot{v}^n(s)$
$dt\times d\loigrandX$-a.s. since $\grandX^u\loigrandX\sim\loigrandX$.

\no An algebraic calculation gives $\loigrandX$-a.s.
$$\rho(-\delta_{\beta_\grandX} v)\circ \grandX^u\rho(-\delta_{\beta_\grandX} u)=1$$
\no Now set $g\in C_b(W,\RR_+)$, we have:
\beaa \EE_{\loigrandX}\left[gL\right]&=&\EE_{\loigrandX}\left[g\circ \grandX^u\right]\\
&=&\EE_{\loigrandX}\left[g\circ
  \grandX^u\rho(-\delta_{\beta_\grandX} v)\circ
  \grandX^u\rho(-\delta_{\beta_\grandX} u)\right]\\
&\leq&\EE_{\loigrandX}\left[g\rho(-\delta_{\beta_\grandX} v) \right]\eeaa

\no So $L\leq\rho(-\delta_{\beta_\grandX}v)$ and since $\EE_{\loigrandX}\left[\rho(-\delta_{\beta_\grandX} v)\right]=1$ we have
$$L\circ \grandX^u\rho(-\delta_{\beta_\grandX} u)=1$$
\no and we can compute $H(\grandX^u\loigrandX|\loigrandX)$:
\beaa H(\grandX^u\loigrandX|\loigrandX)&=&\EE_{\loigrandX}[L\log L]\\ &=&\EE_{\loigrandX}[\log L\circ \grandX^u]\\
&=&\EE_{\loigrandX}\left[-\log\rho(-\delta_{\beta_\grandX} u)\right]\\
&=&\frac{1}{2}\EE_{\loigrandX}[|u_H^2|]\eeaa
\nqed

\section{\bf{Approximation of absolutely continuous measures}}

\begin{theorem}
\label{1repr}
If $\theta\sim \loigrandX$ is such that there exists $r>1$ such that
$$\frac{d\theta}{d\loigrandX}\log\frac{d\theta}{d\loigrandX}\in L^1(\loigrandX)$$
\no and
$$\log\frac{d\theta}{d\loigrandX}\in L^r(\loigrandX)$$
\no there exists $(u_n)\in L^{\infty}_a(\loigrandX,H)^\NN$ such that for
every n,
\beaa &&\grandX^{u_n}\loigrandX\sim \loigrandX\\
&&\grandX^{u_n}\;\; is\;\; \loigrandX-a.s.\;\; invertible\\
&&\frac{d\grandX^{u_n}\loigrandX}{d\loigrandX}\log\frac{d\grandX^{u_n}\loigrandX}{d\loigrandX}\rightarrow \frac{d\theta}{d\loigrandX}\log
\frac{d\theta}{d\loigrandX}\;\; in\;L^1(\loigrandX)\\
&&\frac{d\grandX^{u_n}\loigrandX}{d\loigrandX}\log \frac{d\theta}{d\loigrandX}\rightarrow\frac{d\theta}{d\loigrandX}\log
\frac{d\theta}{d\loigrandX}\;\; in\;L^1(\loigrandX)\eeaa.
\end{theorem}

\nproof Denote
\beaa L&=&\frac{d\theta}{d\loigrandX}\eeaa
\no Eventually sequentializing afterward, we have to prove that for any
$\epsilon>0$, there exists $u\in L^\infty_a(\loigrandX,H)$ such that $\grandX^u\loigrandX\sim \loigrandX$,
$\grandX^u$ is $\loigrandX$-a.s. invertible and
\beaa \EE_{\loigrandX}\left[\left|\frac{d\grandX^u\loigrandX}{d\loigrandX}\log \frac{d\grandX^u\loigrandX}{d\loigrandX}-L\log
    L\right|\right]&\leq&\epsilon\\
 \EE_{\loigrandX}\left[\left|\frac{d\grandX^u\loigrandX}{d\loigrandX}\log L-L\log
    L_1\right|\right]&\leq&\epsilon\eeaa
\no The proof is divided in five steps.\newline
\no Step 1 : We approximate L with a density that is both lower and
upper bounded.\newline
\no Denote
\beaa\phi_n&=&\min(L,n)\\
L_n&=&\frac{\phi_n}{\EE_{\loigrandX}\left[\phi_n\right]}\eeaa
\no The monotone convergence theorem ensures that
$\EE_{\loigrandX}\left[\phi_n\right]\rightarrow 1$ so for any $\alpha\in (0,1)$,
there exists some $n_\alpha\in\NN$ such that for any $n\geq n_\alpha$,
$$\EE_{\loigrandX}\left[\phi_n\right]\geq\alpha$$
\no $(L_n\log L_n)$ converges $\loigrandX$-a.s. to $L\log L$ and if $n\geq n_\alpha
$ and
\beaa
\left|L_n\log
  L_n\right|&=&\left|\frac{\phi_n}{\EE_{\loigrandX}\left[\phi_n\right]}\log
    \frac{\phi_n}{\EE_{\loigrandX}\left[\phi_n\right]}\right|1_{\frac{\phi_n}{\EE_{\loigrandX}\left[\phi_n\right]}\leq 1}+
\left|\frac{\phi_n}{\EE_{\loigrandX}\left[\phi_n\right]}\log
    \frac{\phi_n}{\EE_{\loigrandX}\left[\phi_n\right]}\right|1_{\frac{\phi_n}{\EE_{\loigrandX}\left[\phi_n\right]}>
    1}\\
&\leq& e^{-1}1_{\frac{\phi_n}{\EE_{\loigrandX}\left[\phi_n\right]}\leq 1}+
\left|\frac{L}{\alpha}\log
   \frac{L}{\alpha}\right|1_{\frac{\phi_n}{\EE_{\loigrandX}\left[\phi_n\right]}>
    1}\\
&\leq& e^{-1}+\left|\frac{L}{\alpha}\log
    \frac{L}{\alpha}\right|\eeaa

\no So the Lebesgue theorem ensures that $(L_n\log L_n)$ converge
toward $L\log L$ in $L^1(\loigrandX)$.
\no Similarly, $(L_n\log L)$ converges $\loigrandX$-a.s. to $L\log L$ and
if $n\leq n_\alpha$,
\beaa\left|L_n\log L\right|\leq \left|\frac{L}{\alpha}\log L\right|\eeaa
\no and the Lebesgue theorem ensures that $(L_n\log L)$ converges to
$L_n\log L$ in $L^1(\loigrandX)$, so there exists $n_0\in\NN$ such that
\beaa \EE_{\loigrandX}\left[\left|L_{n_0}\log L_{n_0}-L\log
    L\right|\right]&\leq&\epsilon\\
\EE_{\loigrandX}\left[\left|L_{n_0}\log
    L-L\log L\right|\right]&\leq&\epsilon\eeaa
\no $\left(\frac{L_{n_0}+a}{1+a}\log \frac{L_{n_0}+a}{1+a}\right)$ converges
$\loigrandX$-a.s. to $L_{n_0}\log L_{n_0}$ when a converges to 0. Set $a\in [0,1]$, we have
\beaa \left|\frac{L_{n_0}+a}{1+a}\log \frac{L_{n_0}+a}{1+a}\right|&=&\left|\frac{L_{n_0}+a}{1+a}\log
\frac{L_{n_0}+a}{1+a}\right|1_{L_{n_0}\leq 1}+\left|\frac{L_{n_0}+a}{1+a}\log
\frac{L_{n_0}+a}{1+a}\right|1_{L_{n_0}>1}\\
&\leq& e^{-1}1_{L_{n_0}\leq 1}+\left|L_{n_0}\log
  L_{n_0}\right|1_{L_{n_0}>1}\\
&\leq& e^{-1}+\left|L_{n_0}\log
  L_{n_0}\right|\eeaa
So the Lebesgue theorem ensures that $\left(\frac{L_{n_0}+a}{1+a}\log \frac{L_{n_0}+a}{1+a}\right)$ converges
to $L_{n_0}\log L_{n_0}$ in $L^1(\loigrandX)$. Similarly,
$\left(\frac{L_{n_0}+a}{1+a}\log L\right)$ converges $\loigrandX$-a.s. to
$L_{n_0}\log L$ and
\beaa \left|\frac{L_{n_0}+a}{1+a}\log L\right|&\leq&\left|(L_{n_0}+1)\log L\right|\eeaa
\no and the Lebesgue theorem ensures that $\left(\frac{L_{n_0}+a}{1+a}\log L\right)$ converges to
$L_{n_0}\log L$ in $L^1(\loigrandX)$ and
there exists $a\in [0,1]$ such that
\beaa \EE_{\loigrandX}\left[\left|\frac{L_{n_0}+a}{1+a}\log
  \frac{L_{n_0}+a}{1+a}-L_{n_0}\log L_{n_0}\right|\right]\leq\epsilon\\
\EE_{\loigrandX}\left[\left|\frac{L_{n_0}+a}{1+a}\log
  L-L_{n_0}\log L\right|\right]\leq\epsilon\eeaa
\no $\frac{L_{n_0}+a}{1+a}$ is both lower-bounded and upper-bounded in $L^\infty(\loigrandX)$, we denote these bounds respectively d and D.

\no Also denote
$$M(t)=\EE_{\loigrandX}\left[\left.\frac{L_{n_0}+a}{1+a}\right|\F_t\right]$$
\no We write
$$M(t)=\exp\left(\int_0^t\dot{\alpha}(s)d\beta_\grandX(s)-\frac{1}{2}\int_0^t\left|\dot{\alpha}(s)\right|^2ds\right)$$
\no with $\alpha\in
L^0_a(\loigrandX,H)$.\newline
\no Step 2 : We prove that $\alpha\in L^2(\loigrandX,H)$.\newline
\no Set
$$T_n=\inf\left\{t\in [0,1], \int_0^t\left|\dot{\alpha}(s)\right|^2ds>n\right\}$$
\no $(T_n)$ is a sequence of stopping times which increases
stationarily toward 1. We have, using
$M=1+\int_0^.\dot{\alpha}(s)M(s)d\beta_\grandX(s)$
\beaa \EE_{\loigrandX}\left[\left(M(t\wedge
      T_n)-1\right)^2\right]&=&\EE_{\loigrandX}\left[\int_0^{t\wedge T_n}\left|\dot{\alpha}(s)\right|^2M(s)^2ds\right]\\
&\geq& d^2\EE_{\loigrandX}\left[\int_0^{t\wedge T_n}\left|\dot{\alpha}(s)\right|^2ds\right]\eeaa
\no so
$$\EE_{\loigrandX}\left[\int_0^{t\wedge T_n}\left|\dot{\alpha}(s)\right|^2ds\right]\leq \frac{1}{d^2}\EE_{\loigrandX}\left[\left(M(t\wedge
      T_n)-1\right)^2\right]\leq \frac{2\left(D^2+1\right)}{d^2}$$
\no hence passing to the limit
$$\EE_{\loigrandX}\left[\int_0^{1}\left|\dot{\alpha}(s)\right|^2ds\right]\leq \frac{2\left(D^2+1\right)}{d^2}$$
\no Step 3 : we approximate $\alpha$ with an element of $L^\infty(\loigrandX,H)$.\newline
\no Define
$$\alpha_n:(t,w)\in [0,1]\times \W\mapsto\int_0^t
\dot{\alpha}(s,w)1_{[0,T_n]}(s,w)ds$$

\no and
\beaa
M^n(t)&=&\exp\left(\int_0^t\dot{\alpha^n}(s)d\beta_\grandX(s)-\frac{1}{2}\int_0^t\left|\dot{\alpha^n}(s)\right|^2ds\right)\eeaa
\no and clearly $(M^n(1)\log M^n(1))$ converges $\loigrandX$-a.s. to $M(1)\log M(1)$, $(M^n(1)\log
L)$ converges to $M(1)\log L$ $\loigrandX$-a.s. and
$M^n(1)=\EE_{\loigrandX}\left[\left. M(1)\right|\F_{T_n}\right]$, so $\loigrandX$-a.s.
\beaa \left|M^n(1)\log M^n(1)\right|&\leq &\max \left(e^{-1},\left|D\log
  D\right|\right)\\
\left|M^n(1)\log L\right|&\leq&\left|D\log L\right|\eeaa
\no so the Lebesgue theorem ensures that $(M^n(1)\log M^n(1))$ converges
to $M(1)\log M(1)$ in $L^1(\loigrandX)$ and $(M_1^n\log
L)$ converges to $M_1\log L$ in $L^1(\loigrandX)$ and there exists
$n\in\NN$ such that
\beaa \left|M^{n}(1)\log M^{n}(1)-M(1)\log M(1)\right|&\leq&\epsilon\\
\left|M^{n}(1)\log L-M(1)\log L\right|&\leq&\epsilon\eeaa

\no Step 4 : we approximate $\alpha^n$ with a retarded shift.\newline
\no For $\eta>0$ set
$$\gamma^{\eta}:(t,w)\in[0,1]\times \W\mapsto
\int_0^t\dot{\alpha}^{n}(s-\eta)(w)1_{s>\eta}ds$$
$$N^{\eta}(t)=\exp\left(\int_0^t\dot{\gamma^{\eta}}(s)d\beta_\grandX(s)-\frac{1}{2}\int_0^t\left|\dot{\gamma^{\eta}}(s)\right|^2ds\right)$$
\no $\left(N^\eta(1)\log N^\eta(1)\right)$
converges in probability to $
M^{n}(1)\log M^{n}(1)$.\newline
\no To prove that $\left(N^\eta(1)\log
  N^\eta(1)\right)$ is uniformly integrable, it is sufficient to prove
it is bounded in any $L^p(\loigrandX)$, set $p>1$
\beaa
\EE_{\loigrandX}\left[\left|N^\eta(1)\right|^{p}\right]&=&\EE_{\loigrandX}\left[\exp\left(p\int_0^1\dot{\gamma^\eta}(s)d\beta_\grandX(s)-\frac{p}{2}\int_0^1\left|\dot{\gamma^\eta}(s)\right|^2ds\right)\right]\\
&=&\EE_{\loigrandX}\left[\exp\left(p\int_0^1\dot{\gamma^\eta}(s)d\beta_\grandX(s)-\frac{p^2}{2}\int_0^1\left|\dot{\gamma^\eta}(s)\right|^2ds\right)\exp\left(\frac{p^2-p}{2}\int_0^1\left|\dot{\gamma^\eta}(s)\right|^2ds\right)\right]\\
&\leq&\EE_{\loigrandX}\left[\exp\left(\int_0^1p\dot{\gamma^\eta}(s)d\beta_\grandX(s)-\frac{1}{2}\int_0^1\left|p\dot{\gamma^\eta}(s)\right|^2ds\right)\exp\left(\frac{p^2-p}{2}n\right)\right]\\
&\leq&\exp\left(\frac{p^2-p}{2}n\right)\eeaa
\no so $\left(N^\eta(1)\log N^\eta(1)\right)$ converges to $M^{n}(1)\log
M^{n}(1)$ in $L^1(\loigrandX)$.
\no Furthermore, using H\"older inequality, we have
$$\EE_{\loigrandX}\left[\left|N^\eta(1)\log L-M^{n}(1)\log
    L\right|\right]\leq\left|N^\eta(1)-M^{n}(1)\right|_{L^{r'}(\loigrandX)}\left|\log
  L\right|_{L^r(\loigrandX)}$$
\no where $\frac{1}{r'}+\frac{1}{r}=1$.\newline
\no Consequently there exists $\eta>0$ such that
\beaa \EE_{\loigrandX}\left[\left|N^\eta(1)\log N^\eta(1)-M^{n}(1)\log
    M^{n}(1)\right|\right]&\leq&\epsilon\\
 \EE_{\loigrandX}\left[\left|N^\eta(1)\log L-M^{n}(1)\log
    L\right|\right]&\leq&\epsilon\eeaa
\no using the triangular inequality, we have
\beaa \EE_{\loigrandX}\left[\left|L\log L-N^\eta(1)\log
  N^\eta(1)\right|\right]&\leq&\EE_{\loigrandX}\left[\left|L\log L-L_{n_0}\log
L_{n_0}\right|\right]\\
&&+\EE_{\loigrandX}\left[\left|L_{n_0}\log
L_{n_0}-\frac{L_{n_0}+a}{1+a}\log\frac{L_{n_0}+a}{1+a}\right|\right]\\
&&+\EE_{\loigrandX}\left[\left|\frac{L_{n_0}+a}{1+a}\log\frac{L_{n_0}+a}{1+a}-M^n(1)\log
M^n(1)\right|\right]\\
&&+\EE_{\loigrandX}\left[\left|M^{n}(1)\log M^{n}(1)-N^\eta(1)\log
N^\eta(1)\right|\right]\\
&\leq&4\epsilon\\
\EE_{\loigrandX}\left[\left|L\log L-N^\eta(1)\log
 L\right|\right]&\leq&\EE_{\loigrandX}\left[\left|L\log L-L_{n_0}\log
L\right|\right]\\
&&+\EE_{\loigrandX}\left[\left|L_{n_0}\log
L-\frac{L_{n_0}+a}{1+a}\log L\right|\right]\\
&&+\EE_{\loigrandX}\left[\left|\frac{L_{n_0}+a}{1+a}\log L-M^n(1)\log
L\right|\right]\\
&&+\EE_{\loigrandX}\left[\left|M^{n}(1)\log L-N^\eta(1)\log
L\right|\right]\\
&\leq&4\epsilon\eeaa

\no Step 5 : We prove that $\grandX^{-\gamma^\eta}$ is
$\loigrandX$-a.s. left-invertible and is the solution to our problem.\newline
\no We know $\loigrandX(\grandX=I_\W)=1$ and that there exists a
measurable function $\Phi$ such that $\grandX=\Phi(\beta_\grandX)\;\;\loigrandX$-a.s., so set $A\subset W$,
such that $\loigrandX(A)=1$ and for every $w\in A$, $\grandX(w)=w$ and
$\grandX(w)=\Phi(\beta_\grandX(w))$.\newline

\no Now set $w_1,w_2\in W$ such that
$\grandX^{-\gamma^\eta}(w_1)=\grandX^{-\gamma^\eta}(w_2)$. We have
$$\beta_\grandX(w_1)-\int_0^.\dot{\gamma^\eta}(s)(w_1)ds=\beta_\grandX(w_2)-\int_0^.\dot{\gamma^\eta}(s)(w_2)ds$$
\no For any $s\in[0,\eta]$, $\beta_\grandX(s,w_1)=\beta_\grandX(s,w_2)$,
$\gamma^\eta$ being adapted to filtration $(\F^{\beta_\grandX}_{t-\eta})$,
it implies that for $t\in[0,2\eta]$
$$\int_0^t\dot{\gamma^\eta}(s,w_1)ds=\int_0^t\dot{\gamma^\eta}(s,w_2)ds$$
\no and
$$\beta_\grandX(t,w_1)=\beta_\grandX(t,w_2)$$
\no An easy iteration shows that $\beta_\grandX(w_1)=\beta_\grandX(w_2)$ hence
\beaa w_1&=&\grandX(w_1)\\
&=&\Phi(\beta_\grandX(w_1))\\
&=&\Phi(\beta_\grandX(w_2))\\
&=&\grandX(w_2)\\
&=&w_2\eeaa
\no So $\grandX^{-\gamma^\eta}$ is $\loigrandX$-a.s. injective and so $\loigrandX$-a.s. left-invertible and it is of
the form $\grandX^{v^\eta}$ with $v^\eta\in G_0(\loigrandX,\beta_\grandX)$. We have, for $f\in C_b(\W)$,
\beaa \EE_{\loigrandX}\left[f\circ \grandX^{v^\eta}\right]&=& \EE_{\loigrandX}\left[f\circ \grandX^{v^\eta}\circ \grandX^{-\gamma^\eta}\rho\left(\delta_{\beta_\grandX}\gamma^\eta\right)\right]\\
&=& \EE_{\loigrandX}\left[fN^\eta(1)\right]\eeaa

\no We have
$$\frac{d \grandX^{v^\eta}\loigrandX}{d\loigrandX}=N^\eta(1)$$
\no So $\grandX^{v^\eta}\loigrandX\sim \loigrandX$ and
$$\grandX^{v^\eta}\circ
\grandX^{-\gamma^\eta}= \grandX^{-\gamma^\eta}\circ \grandX^{v^\eta}\;\;
\loigrandX-a.s.$$
\no and
\beaa\left|v^\eta\right|_H^2&=&\left|\gamma^\eta\right|_H^2 \circ
  \grandX^{v^\eta}\\
&\leq&n\eeaa
\no $\loigrandX$-a.s. since $\grandX^{v^\eta}\loigrandX\sim \loigrandX $ and
$\left|\gamma^\eta\right|_H^2\leq n\;\;\loigrandX$-a.s., hence $v^\eta\in L^\infty_a(\loigrandX,H)$.
\nqed

\section{\bf{Variational problem}}
\no As stated in the beginning, we aim to provide a variational
formulation of $-\log\EE_{\loigrandX}\left[e^{-f}\right]$ . This first
result is from \cite{art}:

\begin{theorem}
Set $f:\W\rightarrow \RR$ a measurable function verifying
$$\EE_{\loigrandX}\left[|f|(1+e^{-f})\right]<\infty$$
\no Denote $\mathcal{P}$ the set of probability measures on $\left(\W,\F_1\right)$ which are absolutely continuous with respect to $\loigrandX$, then
$$-\log\EE_{\loigrandX}\left[e^{-f}\right]=\inf_{\theta\in
  \mathcal{P}}\left(\EE_\theta[f]+H(\theta|\loigrandX)\right)$$
\no and the unique supremum is attained at the measure
$$d\theta_0=\frac{e^{-f}}{\EE_{\loigrandX}\left[e^{-f}\right]}d\loigrandX$$
\end{theorem}

\begin{proposition}
\label{1ineqvar}
\no Set $f:\W\rightarrow\RR$ a measurable function
verifying $\EE_{\loigrandX}\left[|f|(1+e^{-f})\right]<\infty$, then
$$-\log\EE_{\loigrandX}\left[e^{-f}\right]\leq\inf_{u\in G_2(\loigrandX,\beta_\grandX)}\EE_{\loigrandX}\left[f\circ \grandX^u+\frac{1}{2}|u|_H^2\right]$$
\end{proposition}

\no Here is the main result.

\begin{theorem}
\label{1varrep}
Set $p>1$ and $f\in L^p(\loigrandX)$ such that
$\EE_{\loigrandX}\left[(|f|+1)e^{-f}\right]<\infty$, then we have
$$-\log\EE_{\loigrandX}\left[e^{-f}\right]=\inf_{u\in L^\infty_a(\loigrandX,H),\grandX^u\;\loigrandX-a.s.\;invertible}\EE_{\loigrandX}\left[f\circ \grandX^u+\frac{1}{2}|u|_H^2\right]$$
\end{theorem}

\nproof Using proposition \ref{1ineqvar}, we have easily
$$-\log\EE_{\loigrandX}\left[e^{-f}\right]\leq\inf_{u\in L^\infty_a(\loigrandX,H),\grandX^u\;\loigrandX-a.s.\;invertible}\EE_{\loigrandX}\left[f\circ \grandX^u+\frac{1}{2}|u|_H^2\right]$$
\no Let $\theta_0$ be the measure on $\W$ defined by
$$d\theta_0=\frac{e^{-f}}{\EE_{\loigrandX}\left[e^{-f}\right]}d\loigrandX$$
\no According to theorem \ref{1repr}, there exists $(u_n)\in
L^\infty_a(\loigrandX,H)^\NN$ such that for every $n\in\NN$, $\grandX^{u_n}$ is
$\loigrandX$-a.s. invertible and that
\beaa \frac{d\grandX^{u_n}\loigrandX}{d\loigrandX}\log \frac{d\grandX^{u_n}\loigrandX}{d\loigrandX} \rightarrow
\frac{d\theta_0}{d\loigrandX}\log \frac{d\theta_0}{d\loigrandX}\\
\frac{d\grandX^{u_n}\loigrandX}{d\loigrandX}\log \frac{d\theta_0}{d\loigrandX} \rightarrow
\frac{d\theta_0}{d\loigrandX}\log \frac{d\theta_0}{d\loigrandX}\eeaa
\no in $L^1(\loigrandX)$.\newline
\no Since $\grandX^{u_n}$
is $\loigrandX$-a.s. invertible, we have
\beaa \EE_{\loigrandX}\left[f\circ \grandX^{u_n}+\frac{1}{2}|u_n|_H^2\right]=
\EE_{\loigrandX}\left[f \frac{d\grandX^{u_n}\loigrandX}{d\loigrandX}\right]+\EE_{\loigrandX}\left[\frac{d\grandX^{u_n}\loigrandX}{d\loigrandX}\log
  \frac{d\grandX^{u_n}\loigrandX}{d\loigrandX}\right]\eeaa
\no When n goes to infinity, we have
$$\EE_{\loigrandX}\left[\frac{d\grandX^{u_n}\loigrandX}{d\loigrandX}\log
  \frac{d\grandX^{u_n}\loigrandX}{d\loigrandX}\right]\rightarrow
\EE_{\loigrandX}\left[\frac{d\theta_0}{d\loigrandX}\log \frac{d\theta_0}{d\loigrandX}\right]$$
\no and since $f=-\log \frac{d\theta_0}{d\loigrandX}-\log\EE_{\loigrandX}\left[e^{-f}\right]$,
\beaa\EE_{\loigrandX}\left[f\frac{d\grandX^{u_n}\loigrandX}{d\loigrandX}\right]&\rightarrow& \EE_{\loigrandX}\left[f\frac{d\theta_0}{d\loigrandX}\right]\eeaa
\no So finally, when n goes to infinity,
\beaa \EE_{\loigrandX}\left[f\circ
  \grandX^{u_n}+\frac{1}{2}|u_n|_H^2\right]&\rightarrow&\EE_{\theta_0}\left[f\right]+H(\theta_0|\loigrandX)\\
&=&-\log\EE_{\loigrandX}\left[e^{-f}\right]\eeaa
\no which conclude the proof.\nqed

\begin{corollary}
Set $f:\W\rightarrow\RR$ such that $\EE_{\loigrandX}\left[(|f|+1)e^{-f}\right]<\infty$ and
$$-\log\EE_{\loigrandX}\left[e^{-f}\right]=\inf_{u\in L^\infty_a(\loigrandX,H),\grandX^u\;\loigrandX-a.s.\;invertible}\EE_{\loigrandX}\left[f\circ \grandX^u+\frac{1}{2}|u|_H^2\right]$$
\no We have
$$-\log\EE_{\loigrandX}\left[e^{-f}\right]=\inf_{u\in G_2(\loigrandX,\beta_\grandX)}\EE_{\loigrandX}\left[f\circ \grandX^u+\frac{1}{2}|u|_H^2\right]$$
\end{corollary}

\begin{theorem}
 Set $f:W\rightarrow\RR$ a measurable function
verifying $\EE_{\loigrandX}\left[|f|(1+e^{-f})\right]<\infty$, then if there
exists some $u\in G_2(\loigrandX,\beta_\grandX)$ such that $\grandX^u$ is
$\loigrandX$-a.s. left-invertible and
$\frac{d\grandX^u\loigrandX}{d\loigrandX}=\frac{e^{-f}}{\EE_{\loigrandX}\left[e^{-f}\right]}$, then
we have
$$-\log\EE_{\loigrandX}\left[e^{-f}\right]=\inf_{u\in G_2(\loigrandX,\beta_\grandX)}\EE_{\loigrandX}\left[f\circ \grandX^u+\frac{1}{2}|u|_H^2\right]$$
\end{theorem}

\nproof Since $\grandX^u$ is $\loigrandX$-a.s. left invertible and that
$\frac{d\grandX^u\loigrandX}{d\loigrandX}=\frac{e^{-f}}{\EE_{\loigrandX}\left[e^{-f}\right]}$. We have
$$\frac{1}{2}\EE_{\loigrandX}\left[|u|_H^2\right]=H(\grandX^u\loigrandX|\loigrandX)=\EE_{\loigrandX}\left[\frac{e^{-f}}{\EE_{\loigrandX}\left[e^{-f}\right]}\log
  \left(\frac{e^{-f}}{\EE_{\loigrandX}\left[e^{-f}\right]}\right)\right]$$
\no and
\beaa\EE_{\loigrandX}\left[f\circ
  \grandX^u+\frac{1}{2}|u|_H^2\right]&=&\EE_{\loigrandX}\left[\frac{e^{-f}}{\EE_{\loigrandX}\left[e^{-f}\right]}f+\frac{e^{-f}}{\EE_{\loigrandX}\left[e^{-f}\right]}\log\left(\frac{e^{-f}}{\EE_{\loigrandX}\left[e^{-f}\right]}\right)\right]\\
&=&-\log\EE_{\loigrandX}\left[e^{-f}\right]\eeaa
\no and we conclude the proof with last proposition.\nqed

\begin{theorem}
Set $f:W\rightarrow \RR$  a measurable function such that
$$-\log\EE_{\loigrandX}\left[e^{-f}\right] = \inf_{u\in
  G_2(\loigrandX,\beta_\grandX)}\EE_\nu\left[f\circ \grandX^u+\frac{1}{2}|u|_H^2\right]$$
\no Denote this infimum $J_*$. It is attained at $u\in G_2(\loigrandX,\beta_\grandX)$ if and only
if $\grandX^u$ is $\loigrandX$-a.s. left-invertible and
$\frac{d\grandX^u\loigrandX}{d\loigrandX}=\frac{e^{-f}}{\EE_{\loigrandX}\left[e^{-f}\right]}$.\newline
\end{theorem}

\nproof The direct implication is given by last theorem.
\no Conversely, if $\grandX^u$ is not $\loigrandX$-a.s. left-invertible,
$H(\grandX^u\loigrandX|\loigrandX)<\frac{1}{2}\EE_{\loigrandX}\left[|u|_H^2\right]$ and
\beaa -\log\EE_{\loigrandX}\left[e^{-f}\right]=\inf_{\theta\in
  \mathcal{P}(W)}\left(\EE_\theta[f]+H(\theta|\loigrandX)\right)&\leq& \inf_{\alpha\in
  G_2(\loigrandX,\beta_\grandX)}\EE_{\grandX^\alpha \loigrandX}\left[f\right]+H(\grandX^\alpha \loigrandX|\loigrandX)\\
&\leq&\EE_{\grandX^u \loigrandX}\left[f\right]+H(\grandX^u \loigrandX|\loigrandX)\\
&<&\EE_{\loigrandX}\left[f\circ \grandX^u+\frac{1}{2}|u|_H^2\right]\eeaa
\no which is a contradiction.\newline
\no We get $\frac{d\grandX^u\loigrandX}{d\loigrandX}=L$ by uniqueness of the minimizing
measure of $\inf_{\theta\in
  \mathcal{P}(W)}\left(\EE_\theta[f]+H(\theta|\loigrandX)\right)$.\newline
\nqed

\vspace{2cm}
\footnotesize{
\noindent
K\'evin HARTMANN, Institut Telecom, Telecom ParisTech, LTCI CNRS D\'ept. Infres, \\
23 avenue d'Italie, 75013, Paris, France\\
kevin.hartmann@polytechnique.org}

\end{document}